\documentclass[twocolumn,twoside]{IEEEtran}
\usepackage{amsmath,amssymb,amsthm,wasysym,epsfig,dsfont,color,subfigure,empheq,graphicx}
\usepackage{enumerate,stfloats,url,algpseudocode,algorithm}
\usepackage[table]{xcolor}

\DeclareMathOperator{\sign}{sign}

\newtheorem{proposition}{Proposition}

\theoremstyle{remark}\newtheorem{remark}{Remark}

\begin{document}
\title{Stochastic Reactive Power Management\\
in Microgrids with Renewables}

\author{
	Vassilis Kekatos,~\IEEEmembership{Member,~IEEE,}
	Gang Wang,~\IEEEmembership{Student Member,~IEEE,}\\
	Antonio J. Conejo,~\IEEEmembership{Fellow,~IEEE,}  and
	Georgios B. Giannakis,~\IEEEmembership{Fellow,~IEEE}
			



}


\maketitle

\begin{abstract}
Distribution microgrids are being challenged by reverse power flows and voltage fluctuations due to renewable generation, demand response, and electric vehicles. Advances in photovoltaic (PV) inverters offer new opportunities for reactive power management provided PV owners have the right investment incentives. In this context, reactive power compensation is considered here as an ancillary service. Accounting for the increasing time-variability of distributed generation and demand, a stochastic reactive power compensation scheme is developed. Given uncertain active power injections, an online reactive control scheme is devised. This scheme is distribution-free and relies solely on power injection data. Reactive injections are updated using the Lagrange multipliers of a second-order cone program. Numerical tests on an industrial 47-bus microgrid and the residential IEEE 123-bus feeder corroborate the reactive power management efficiency of the novel stochastic scheme over its deterministic alternative, as well as its capability to track variations in solar generation and household demand.
\end{abstract}


\section{Introduction}\label{sec:intro}
Medium- and low-voltage power grids nowadays are undergoing a transformative change to microgrids. Renewable generation and elastic loads are uncertain, power flows are frequently reversed, and bus voltage magnitudes can fluctuate considerably. For example, the power generated by a photovoltaic (PV) network with intermittent cloud coverage can vary by 15\% of its nameplate capacity within one-minute intervals~\cite{Turitsyn11}. Different from transmission grids, bus voltage magnitudes in distribution grids are markedly affected by active power variations. On a clear day, solar generation may easily exceed local demand (especially at midday off-peak hours) and cause over-voltages~\cite{Carvalho}; whereas overnight vehicle charging could lead to serious voltage sags~\cite{Guille}.

Given active power injections, reactive power management aims at controlling reactive injections so that power losses over distribution lines are minimized while bus voltage magnitudes are maintained within the prescribed limits, e.g., $\pm 5\%$ of their nominal values. Traditionally, reactive power management is achieved via tap-changing under load (TCUL) transformers, step voltage regulators (SVR), shunt capacitors and reactors, and static var compensators (SVC)~\cite{Kersting}; see for example \cite{BaldickWu90,Genetic} for related control algorithms. Operational costs, discrete control actions, and slow response times are the factors limiting the use of such devices alone for voltage regulation in distribution systems with renewables~\cite{FCL}. Reactive power management becomes even more challenging in microgrids operating in islanded mode due to the lack of centralized fast-reacting generators~\cite{Microgrids}. On the other hand, subsidizing reactive power control by distributed generation (DG) units has been advocated as a viable solution~\cite{Rogers10,Turitsyn11}. 

Although prohibited by current standards~\cite{IEEE1547}, the power electronics of PVs can be commanded to provide reactive injections as well; see \cite{Robbins} and references therein. For this reason, reactive power compensation via DG units has been an active research area lately. A multi-agent approach is proposed in~\cite{Markabi}, while voltage regulation is cast as a learning problem in~\cite{Vaccaro}. Control policies based on approximate models are developed in~\cite{Turitsyn11}; and a successive convex approximation is adopted in~\cite{Deshmukh} for voltage regulation. Upon linearizing the power flow equations, a two-layer decentralized scheme is proposed in~\cite{Robbins}. Another decentralized consensus-type algorithm is pursued in~\cite{Saverio} after approximating power losses as a quadratic function of reactive power injections. Localized (re)active injection updates are reported in~\cite{FarivarCDC},~\cite{GermanCode},~\cite{ZDGT13}.

All previous schemes build on approximate grid models. Being an instance of the optimal power flow (OPF) problem, reactive power management is a non-convex problem, yet several convex relaxations have been proposed~\cite{ Low14}. In radial distribution grids, OPF can be surrogated by a semidefinite program (SDP)~\cite{Bai08,Lavaei12}; or by a second-order cone program (SOCP) using either polar coordinates~\cite{Jabr06}, or the branch flow model~\cite{BW1,BW2,FL1}. A one-to-one mapping between their feasible sets proves the equivalence of the two relaxations~\cite{BoseAllerton}, and advocates using the SOCP one due to its simplicity. Sufficient conditions guaranteeing the exactness of the convex relaxation (i.e., that solving the relaxed problem is equivalent to solving the original non-convex one) have been developed; see  \cite{Low14} for a review. Regarding reactive power compensation, a distributed algorithm based on the SDP relaxation has been developed in~\cite{ZLGT13}, and a centralized approach for inverter VAR control using the SOCP relaxation has been devised in~\cite{FCL}.

The approaches so far assume that active power injections are precisely known and remain unchanged throughout the reactive control period. However, such assumptions are less realistic in future microgrids with high penetration of renewables. Our first contribution is a stochastic framework for reactive power management. We consider a radial microgrid where several DG units with reactive power control capabilities have been integrated. The grid operation is divided into short time intervals. At every interval, a microgrid controller collects active nodal injections and decides the reactive power to be injected by controllable DG units. (Re)active load demands and renewable generation are known only via noisy and delayed estimates, and are hence, modeled as stochastic processes. Different from the power loss minimization in \cite{NAPS2014}, an ancillary voltage regulation market is formulated here: PV owners are reimbursed for providing reactive power support. Reactive injections from PVs are set as the minimizers of an \emph{expected} reactive power compensation cost. 

As a second contribution, the derived optimization problem is solved using a provably convergent stochastic approximation algorithm. It is further shown that a subgradient of the involved cost is computed via the dual SOCP problem and reactive PV injections are updated by a simple thresholding rule. Numerical results on industrial and residential microgrids with real solar generation and demand data corroborate the efficacy of the novel reactive power management scheme.

The rest of the paper is outlined as follows: After the branch flow model is presented in Section~\ref{sec:model}, the problem of stochastic reactive power compensation is formulated in Section~\ref{sec:problem}. A stochastic approximation algorithm is developed in Section~\ref{sec:solver}, its performance advantage over an instantaneous reactive control scheme is supported by numerical tests in Section~\ref{sec:simulations}, and conclusions are drawn in Section~\ref{sec:conclusions}.

Regarding notation, lower- (upper-) case boldface letters denote column vectors (matrices), with the exception of line power flow vectors $(\mathbf{P},\mathbf{Q})$. Calligraphic symbols are reserved for sets. Prime stands for vector and matrix transposition. Vectors $\mathbf{0}$ and $\mathbf{e}_n$ denote the all-zeros and the $n$-th canonical vector, respectively. Symbol $\|\mathbf{x}\|_2$ denotes the $\ell_2$-norm of $\mathbf{x}$.
 
 \begin{figure}[t]
\centering
\includegraphics[scale=0.3]{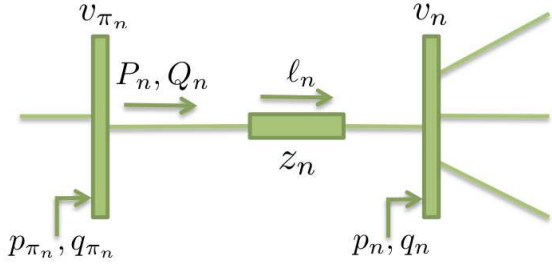}
\caption{Bus $n$ is connected to its unique parent $\pi_n$ via line $n$.}
\label{fig:diagram}
\end{figure}
 
\section{System Model}\label{sec:model}
Consider a microgrid consisting of $N+1$ buses. For operational and architectural simplicity, the microgrid is assumed to be radial and it can thus be modeled by a tree graph $\mathcal{T}:=(\mathcal{N}_o,\mathcal{L})$, where $\mathcal{N}_o:=\{0,1,\ldots,N\}$ denotes the set of nodes (buses), and $|\mathcal{L}|=N$ is the cardinality of the edge set $\mathcal{L}$. The tree is rooted at the substation bus indexed by $n=0$. For every bus $n\in\mathcal{N}_o$, let $v_n$ be the squared voltage magnitude at bus $n$ and $p_n+jq_n$ the complex power injected into bus $n$. Notice that every non-root bus $n\in\mathcal{N}=\{1,\ldots,N\}$ has a unique parent bus that will be denoted by $\pi_n$. Hence, the directed edge $(\pi_n,n)\in\mathcal{L}$ corresponding to the distribution line feeding bus $n$ will be simply indexed by $n$; see Fig.~\ref{fig:diagram}. Let also $z_n=r_n+jx_n$ and $\ell_n$ denote the line impedance and the squared current magnitude on line $n$, respectively. If $P_n+jQ_n$ is as the complex power flow on line $n$ seen at the sending end $\pi_n$, the so termed \emph{branch flow model} is described by the equations~\cite{BW1,BW2},
\begin{align}
p_n&=\sum_{k\in\mathcal{C}_n}P_k  - (P_n -r_n \ell_n)\label{eq:mp}\\
q_n&=\sum_{k\in\mathcal{C}_n}Q_k  - (Q_n -x_n \ell_n)\label{eq:mq}\\
v_n&=v_{\pi_n}+(r_n^2+x_n^2)\ell_n - 2(r_nP_n+x_nQ_n)\label{eq:mv}\\
\ell_n&=\frac{P_n^2+Q_n^2}{v_{\pi_n}}\label{eq:ml}
\end{align}
for all $n\in \mathcal{N}$, where $\mathcal{C}_n:=\{k\in\mathcal{N}:\pi_k=n\}$ is the set of the children nodes for bus $n$. Equations \eqref{eq:mp}-\eqref{eq:mq} follow from power conservation; \eqref{eq:mv} is derived upon squaring second Kirchoff's law; and \eqref{eq:ml} from current computations. The branch flow model is essentially derived from the full AC model, after eliminating voltage and current phases. The model is accompanied with the initial conditions $v_0=1$, $p_0=\sum_{k\in\mathcal{C}_0}P_k$, and $q_0=\sum_{k \in \mathcal{C}_0} Q_k$.

The active and reactive power injection at bus $n$ can be decomposed into its generation and consumption components as $p_n=p_n^g-p_n^c$ and $q_n=q_n^g-q_n^c$. For a purely load bus, there is no generation $(p_n^g=q_n^g=0)$, the consumed active power is $p_n^c\geq 0$, and its reactive power $q_n^c\geq 0$ is typically related to $p_n^c$ via a constant power factor. A DG bus (e.g., an industrial facility equipped with rooftop solar panels or a wind turbine) not only consumes power denoted by $p_n^c$ and $q_n^c$, but it can also generate active power $p_n^g\geq 0$, and provide reactive support $q_n^g$ which can be positive or negative. For a bus hosting a shunt capacitor only, $p_n=q_n^c=0$ and $q_n^g>0$.

For notational simplicity, let us collect all nodal quantities related to non-root buses in vectors $\mathbf{p}:=[p_1~\cdots~p_N]'$, $\mathbf{q}:=[q_1~\cdots~q_N]'$, and $\mathbf{v}:=[v_1~\cdots~v_N]'$. Likewise, for line quantities define vectors $\mathbf{P}:=[P_1~\cdots~P_N]'$, $\mathbf{Q}:=[Q_1~\cdots~Q_N]'$, and $\boldsymbol{\ell}:=[\ell_1~\cdots~\ell_N]'$. Bus voltage magnitudes are allowed to lie within a prespecified range (typically a $\pm 5\%$ of their nominal value), yielding the voltage limits $v_n\in[\underline{v}_n,\overline{v}_n]$ for all $n\in \mathcal{N}$. Upon setting $\underline{\mathbf{v}}:=[\underline{v}_1~\cdots~\underline{v}_N]'$ and $\overline{\mathbf{v}}:=[\overline{v}_1~\cdots~\overline{v}_N]'$, voltage regulation constraints can be compactly expessed as
\begin{equation}\label{eq:Vset}
\mathbf{v}\in\mathcal{V}:=\{\mathbf{v}:\underline{\mathbf{v}}\leq \mathbf{v} \leq\overline{\mathbf{v}}\}.
\end{equation}
Building on \eqref{eq:mp}-\eqref{eq:Vset}, our stochastic reactive control scheme is formulated next.

\section{Problem Formulation}\label{sec:problem}
In the envisioned microgrid operation scenario, active power is managed at a coarse timescale. For example, a power dispatch is issued for the next 24 hours through a day-ahead market. Active power adjustments are implemented on a 5- or 10-minute basis via a real-time market. Both the hourly and the real-time market active power dispatch could depend on the cost of dispatchable generators and predictions on renewable energy within the microgrid, as well as the costs of power exchanges with a main grid. Together with this hourly active power schedule, the microgrid controller manages reactive power by controlling transformers, shunt capacitors, SVRs, and SVCs~\cite{Kersting,BaldickWu90}. Nonetheless, slow response times and switching limitations render such devices inadequate for very fast reactive power control. The power electronic interfaces found in DG units, such as PV inverters, provide a viable solution for near real-time reactive power management~\cite{Turitsyn11},~\cite{ZLGT13}.

Reactive power compensation occurs over time intervals indexed by $t$. These intervals could either coincide with real-time market periods (e.g., 5 minutes), or be even shorter (30 seconds), depending on the variability of active powers and cyber resources (sensing, communication, and computation delays). If $(\mathbf{p}_t,\mathbf{q}_t)$ are the active and reactive power injections in all but the root buses during control period $t$, the power loss on distribution lines is expressed as
\begin{equation}\label{eq:losses}
f(\mathbf{p}_t,\mathbf{q}_t)=\sum_{n=0}^N p_{n,t}=\sum_{n=1}^N r_n \ell_{n,t}
\end{equation}
where the second equality follows from \eqref{eq:mp}. Recalling that $\mathbf{q}_t:=\mathbf{q}_t^g-\mathbf{q}_t^c$, define for notational brevity
\begin{equation}\label{eq:ft}
f_t(\mathbf{q}^g):=f(\mathbf{p}_t,\mathbf{q}^g-\mathbf{q}_t^c).
\end{equation}

Given active injections $\mathbf{p}_t$ and reactive demands $\mathbf{q}_t^c$, conventional reactive power management aims at choosing $\mathbf{q}_t^g$ so that power losses are minimized and voltages are maintained within $\mathcal{V}$. Concretely, reactive power management could be stated as finding
\begin{equation}\label{eq:instantaneous}
\tilde{\mathbf{q}}^g_t:=\arg\min_{\mathbf{q}^g\in \mathcal{Q}}~f_t(\mathbf{q}^g)
\end{equation}
where $\mathcal{Q}$ is the reactive feasible region to be delineated later. Injecting $\tilde{\mathbf{q}}^g_t$ at time $t+1$ would be the optimal control action under two operational conditions:\\
\hspace*{1em} \textbf{(C1)} $(\mathbf{p}_t,\mathbf{q}_t^c)$ are precisely known, and\\
\hspace*{1em} \textbf{(C2)} they remain constant throughout period $t+1$.\\
Yet such conditions are hardly met in microgrids: renewable DG entails time-varying active and reactive injections. In low-inertia microgrids, the lack of droop controllers challenges further voltage regulation. It is worth noting that even if $(\mathbf{p}_t,\mathbf{q}_t^c)$ are relatively constant over periods $t$ and $t+1$, the microgrid controller has only their noise-contaminated observations (direct measurements or delayed state estimates).

To overcome these difficulties, a stochastic optimization approach is pursued here. The active and reactive power injections $(\mathbf{p}_t,\mathbf{q}_t^c)$ realized over an hour or over a real-time market interval are modeled as \emph{stochastic processes} drawn independently across time from a probability density function (pdf): Injections $\{\mathbf{p}_t\}$ could be modeled as the sum of a nominal $\mathbf{p}^o$ and deviations $\{\boldsymbol{\epsilon}_t\}$ that are assumed independent over time; and likewise for $\{\mathbf{q}_t^c\}$. A meaningful stochastic control scheme could entail minimizing the \emph{average} power loss as~\cite{NAPS2014}
\begin{equation}\label{eq:ensemble}
\hat{\mathbf{q}}^g:=\arg\min_{\mathbf{q}^g\in \mathcal{Q}}~
\mathbb{E}[f_t(\mathbf{q}^g)]
\end{equation}  
where the expectation is over time $t$, or more precisely, over $(\mathbf{p}_t,\mathbf{q}_t^c)$. Rather than implementing the unreliable and possibly obsolete instantaneous decisions $\tilde{\mathbf{q}}_t^g$ of \eqref{eq:instantaneous}, problem \eqref{eq:ensemble} is expected to yield smoother control actions. Distinct from \cite{NAPS2014} where PVs were providing voltage regulation at no charge, reactive compensation is interpreted here as an ancillary service. Before elaborating this service, the injection region $\mathcal{Q}$ should be understood first.

\emph{Reactive power injection region:} Choosing $\mathcal{Q}$ requires understanding the reactive control capabilities of PVs~\cite{Turitsyn11}. Consider a solar panel located at bus $n$ with nameplate active power capacity $\overline{p}_{n}^g$, and its inverter having apparent power capability $s_n$. Because PVs are currently restricted to operate at unity power factor~\cite{IEEE1547}, their inverters are typically designed so that $s_n=\overline{p}_n^g$. If $p_{n,t}^g$ is the PV output at time $t$, the inverter could compensate $q_{n,t}^g$ constrained as $|q_{n,t}^g|\leq \sqrt{s_n^2-(p_{n,t}^g)^2}$. This design constraint introduces two practical concerns: First, the reactive injection region becomes time-varying thus complicating \eqref{eq:ensemble}. Second, when $p_{n,t}^g=\overline{p}_{n}^g$ (at maximum solar output), no reactive power can be provided although at those instances it may be needed. 

For these reasons, PV inverters have been advocated to be oversized over their panel nameplate capacity so that $s_n>\overline{p}_n^g$; c.f.~\cite{Turitsyn11}. By choosing for example $s_n=1.1\overline{p}_n^g$  and limiting reactive power compensation to $\sqrt{s_n^2-(\overline{p}_n^g)^2}$ rather than $\sqrt{s_n^2-(p_{n,t}^g)^2}$, the inverter $n$ can provide reactive power support with $|q_{n,t}^g|\leq 0.45\overline{p}_n^g$, \emph{regardless} of the instantaneous PV output $p_{n,t}^g$. Under this policy, the reactive injection region $\mathcal{Q}$ is the time-invariant convex set
\begin{equation}\label{eq:Q}
\mathcal{Q}:=\left[\underline{\mathbf{q}}^g,\overline{\mathbf{q}}^g\right]
\end{equation}
where $\underline{\mathbf{q}}^g:=[\underline{q}_1^g~\cdots~\underline{q}_N^g]'$ and $\overline{\mathbf{q}}^g:=[\overline{q}_1^g~\cdots~\overline{q}_N^g]'$. If $\mathcal{N}_q\subseteq \mathcal{N}$ is the subset of buses with controllable reactive injections, then $\overline{q}_n^g=-\underline{q}_n^g\geq 0$ for $n\in\mathcal{N}_q$, and $\overline{q}_n^g=\underline{q}_n^g=0$ for $n\notin\mathcal{N}_q$.

Although the aforementioned scheme could be technologically feasible, PV owners have to invest on oversized inverters. As a financial incentive, PV sites with reactive power compensation capabilities can participate in an ancillary voltage regulation market and be reimbursed for their reactive power support. Specifically, let $\tilde{c}_n\geq 0$ [in $\cent$/kVar \& h] be the price for reactive power support at bus $n\in\mathcal{N}_q$; and $\tilde{c}_0> 0$ [in $\cent$/kWh] the price at which the microgrid buys (or sells) active power $p_{0,t}$ from (to) the main grid. If prices are constant throughout the real-time market interval or longer periods, a market-based extension of \eqref{eq:ensemble} could be formulated as
\begin{equation}\label{eq:market_tilde}
\hat{\mathbf{q}}^g:=\arg\min_{\mathbf{q}^g\in \mathcal{Q}}~\tilde{c}_0 \mathbb{E}[f_t(\mathbf{q}^g)] + \sum_{n\in\mathcal{N}_q} \tilde{c}_n |q_n^g|
\end{equation}
where the microgrid controller trades power losses for reactive power support by PVs. The ancillary market in \eqref{eq:market_tilde} can be equivalently expressed as
\begin{equation}\label{eq:market}
\hat{\mathbf{q}}^g:=\arg\min_{\mathbf{q}^g\in \mathcal{Q}} ~\mathbb{E}[f_t(\mathbf{q}^g)] + \sum_{n=1}^N c_n |q_n^g|
\end{equation}
where $c_n:=\tilde{c}_n/\tilde{c}_0$ are the normalized prices for $n\in\mathcal{N}_q$, and $c_n:=0$ for $n\notin \mathcal{N}_q$. Even if the joint pdf of $(\mathbf{p}_t,\mathbf{q}_t^c)$ were known, evaluating the expectation in \eqref{eq:ensemble} or \eqref{eq:market} would be non-trivial. To practically solve these two policies, a stochastic approximation approach is pursued next. The focus will be on solving \eqref{eq:market}, yet that is without loss of generality since \eqref{eq:ensemble} is the special case of \eqref{eq:market} where $\{c_n=0\}_{n\in\mathcal{N}}$.

\section{Stochastic Approximation Solver}\label{sec:solver}
Leveraging recent advances in online convex optimization (see for instance \cite{COMID}), the problem in \eqref{eq:market} can be cast in a stochastic approximation framework. Successive estimates $\{\hat{\mathbf{q}}_t^g\}_{t}$ for the minimizer $\hat{\mathbf{q}}^g$ are iteratively found as soon as a new datum $(\mathbf{p}_t,\mathbf{q}_t^c)$ becomes available. Specifically, adopting the composite objective mirror descent approach of~\cite{COMID}, reactive injection iterates $\hat{\mathbf{q}}_t^g$'s are updated as the minimizers of the convex problem
\begin{equation}\label{eq:sa}
\hat{\mathbf{q}}_t^g :=\arg\min_{\mathbf{q}^g\in\mathcal{Q}}~\mathbf{g}_t'\mathbf{q}^g + \tfrac{1}{2\eta_t} \|\mathbf{q}^g-\hat{\mathbf{q}}_{t-1}^g\|_2^2 + \sum_{n=1}^N c_n |q_n^g|
\end{equation}
where $\mathbf{g}_t$ is an arbitrary subgradient of $f_t(\mathbf{q}^g)$ evaluated at $\hat{\mathbf{q}}_{t-1}^g$ and $\eta_t>0$ is an appropriately selected step size. Recall that the subgradient generalizes the notion of gradient to non-differentiable functions. For a convex function $f:\mathbb{R}^n\rightarrow \mathbb{R}$, any vector $\mathbf{g}\in\mathbb{R}^n$ satisfying the inequality $f(\mathbf{y})\geq f(\mathbf{x}) + \mathbf{g}'(\mathbf{x}-\mathbf{y})$ for every $\mathbf{y}$ in the domain of $f$, constitutes a \emph{subgradient} of $f$ at $\mathbf{x}$. The set of all subgradients is termed the \emph{subdifferential} of $f$ at $\mathbf{x}$, and is denoted by $\partial f (\mathbf{x})$; see also \cite[Sec.~2.5]{Ruszczynski}.

The rationale behind stochastic approximation algorithms is to first surrogate the original cost in \eqref{eq:market} by its instantaneous approximation to yield
\begin{equation}\label{eq:instant}
\min_{\mathbf{q}^g\in\mathcal{Q}}~f_t(\mathbf{q}^g) + \sum_{n=1}^N c_n|q_n^g|.
\end{equation}
Notice that the minimizer of \eqref{eq:instant} is the per-time optimal reactive injection. But instead of solving \eqref{eq:instant}, online optimization algorithms minimize a locally tight upper bound of the cost in \eqref{eq:instant}. Such a bound can be obtained by maintaining $\sum_{n=1}^N c_n|q_n^g|$ and linearizing $f_t(\mathbf{q}^g)$ at the previous iterate as $f_t(\hat{\mathbf{q}}_{t-1}^g)+\mathbf{g}_t'(\mathbf{q}-\hat{\mathbf{q}}_{t-1}^g) + \tfrac{1}{2\eta_t}\|\mathbf{q}-\hat{\mathbf{q}}_{t-1}^g\|_2^2$ for a proper $\eta_t>0$. After ignoring constant terms, the update in \eqref{eq:sa} follows.

To practically implement the stochastic reactive control scheme of \eqref{eq:sa}, two issues need to be resolved: finding the minimizer of \eqref{eq:sa} and calculating the subgradient $\mathbf{g}_t$ involved.

\subsection{Closed-Form Minimizer for \eqref{eq:sa}}\label{subsec:closed}
Suppose a subgradient $\mathbf{g}_t$ has been found. Upon completing the square, the optimization in \eqref{eq:sa} can be written as
\begin{equation}\label{eq:sa2}
\hat{\mathbf{q}}_t^g :=\arg\min_{\mathbf{q}^g\in\mathcal{Q}} ~
\frac{1}{2} \|\mathbf{q}^g-\mathbf{y}_t\|_2^2 + \eta_t \sum_{n=1}^N c_n |q_n^g|
\end{equation}
where $\mathbf{y}_t:=\hat{\mathbf{q}}_{t-1}^g -\eta_t \mathbf{g}_t$. Note that solving \eqref{eq:sa2} decouples over the entries of $\mathbf{q}^g$ as 
\begin{equation}\label{eq:sa3}
\hat{q}_{n,t}^g :=\arg\min_{\underline{q}_{n}^g\leq q_{n}^g \leq \overline{q}_{n}^g} ~
\frac{1}{2} \left(q_n^g-y_{n,t}\right)^2 + \eta_t c_n |q_n^g|
\end{equation}
where $y_{n,t}$ is the $n$-th entry of $\mathbf{y}_t$. By using the Karush-Kuhn-Tucker conditions for the univariate minimization in \eqref{eq:sa3}, the following result is shown in the Appendix.

\begin{proposition}\label{pro:closed}
The minimizer of \eqref{eq:sa3} is expressed in closed form as
\begin{equation}\label{eq:closed}
\hat{q}_{n,t}^g =\left\{ \begin{array}{ll}
\overline{q}_{n}^g,		&	y_{n,t}> \overline{q}_{n}^g + \eta_t c_n\\
 y_{n,t} - \eta_t c_n,		&	\eta_t c_n < y_{n,t} \leq \overline{q}_{n}^g + \eta_t c_n\\
0,						&	-\eta_t c_n \leq y_{n,t} \leq \eta_t c_n\\
 y_{n,t} + \eta_t c_n,		&	 \underline{q}_{n}^g - \eta_t c_n \leq y_{n,t} < -\eta_t c_n\\
\underline{q}_{n}^g,		&	y_{n,t}< \underline{q}_{n}^g - \eta_t c_n
\end{array}\right..
\end{equation}
\end{proposition}

The rule of \eqref{eq:closed} implies that if $|y_{n,t}|$ is smaller than $\eta_t c_n$, there is no reactive injection into bus $n$. When $|y_{n,t}|$ is large, its reactive injection saturates. Otherwise, $\hat{q}_{n,t}^g=y_{n,t}-\sign(y_{n,t})\eta_t c_n$. Therefore, once a $\mathbf{g}_t$ belonging to the subdifferential $\partial f_t(\hat{\mathbf{q}}_{t-1}^g)$ has been found, $\hat{\mathbf{q}}_t^g$ can be easily obtained from \eqref{eq:closed}.

\subsection{Efficient Subgradient Computation}\label{subsec:subgradient}
Before finding a subgradient $\mathbf{g}_t$, an alternative representation for $f_t(\mathbf{q}^g)$ is derived first. Recall that $f_t(\mathbf{q}^g)$ is the power loss on distribution lines when injecting $(\mathbf{p}_t,\mathbf{q}^g-\mathbf{q}_t^c)$ into the distribution grid. Provided that $(\mathbf{p}_t,\mathbf{q}^g-\mathbf{q}_t^c)$ is feasible, function $f_t(\mathbf{q}^g)$ depends on the underlying grid operating point $(p_0,q_0,\mathbf{p},\mathbf{q},\mathbf{P},\mathbf{Q},\boldsymbol{\ell},\mathbf{v})$ [cf.~\eqref{eq:losses}]. Finding this point requires solving the nonlinear equations in \eqref{eq:mp}-\eqref{eq:ml}, while guaranteeing that voltages are maintained in the desired range $(\mathbf{v}\in \mathcal{V})$. Solving this set of nonlinear equations and linear inequalities is non-trivial. Under practical operating conditions, the solution has been claimed to be unique~\cite{Baran3}. Even if multiple solutions exist, the grid operating point attaining the smallest loss for the same $(\mathbf{p}_t,\mathbf{q}^g-\mathbf{q}_t^c)$ can be found as described next.

If the equalities in \eqref{eq:ml} are relaxed to inequalities, then $(p_0,q_0,\mathbf{P},\mathbf{Q},\boldsymbol{\ell},\mathbf{v})$ lies in a convex set; see e.g.,~\cite{FL1,FCL}. This convex set is represented by the linear equalities \eqref{eq:mp}-\eqref{eq:mv}, the set $\mathcal{V}$, and the second-order cone constraints $\left\{P_n^2+Q_n^2\leq \ell_n v_{\pi_n}\right\}_{n\in\mathcal{L}}$. Under different technical conditions (see~\cite{Low14} and references therein), the minimizer of the convex problem
\begin{subequations}\label{eq:socp}
\begin{align}
f(\mathbf{p},\mathbf{q})=\min_{{\mathbf{P},\mathbf{Q} \atop \boldsymbol{\ell}, \mathbf{v}}} &~\sum_{n=1}^L r_n \ell_n\label{eq:scost}\\ 
\textrm{s.to}~&p_n=\sum_{k\in\mathcal{C}_n}P_k  - (P_n -r_n \ell_n),~n\in\mathcal{N}\label{eq:sp}\\
&q_n=\sum_{k\in\mathcal{C}_n}Q_k  - (Q_n -x_n \ell_n),~n\in\mathcal{N}\label{eq:sq}\\
&v_n=v_{\pi_n}+(r_n^2+x_n^2)\ell_n- 2(r_nP_n+x_nQ_n),\nonumber\\
&\vspace*{14em}n\in \mathcal{N}\label{eq:sv}\\
&\ell_n\geq \frac{P_n^2+Q_n^2}{v_{\pi_n}},~n\in\mathcal{L}\label{eq:sl}\\
&\mathbf{v}\in \mathcal{V}\label{eq:sV}
\end{align}
\end{subequations}
satisfies the SOCP constraints in \eqref{eq:sl} with equality. When this occurs, the convex relaxation is said to be \emph{exact}. To summarize, when the relaxation is exact, the optimum value of \eqref{eq:socp} equals the loss experienced under injections $(\mathbf{p},\mathbf{q})$.

Henceforth, the following assumptions will be adopted:\\
\hspace*{1em}\textbf{(A1)} \textit{The convex relaxation in \eqref{eq:socp} is exact}.\\
\hspace*{1em}\textbf{(A2)} \textit{There exists a feasible $(\mathbf{P},\mathbf{Q},\boldsymbol{\ell},\mathbf{v})$ for \eqref{eq:socp} satisfying constraints \eqref{eq:sl} with strict inequality.}

Albeit assumptions (A1)-(A2) are not supported analytically here, they are verified throughout our numerical tests. All instances of \eqref{eq:socp} encountered in Section~\ref{sec:simulations} were exact. Additionally, when for these instances the cost in \eqref{eq:scost} was maximized rather than minimized, the resultant maximizers satisfied \eqref{eq:sl} with strict inequality; thus, numerically verifying (A2).

Under (A1), the $(\mathbf{P},\mathbf{Q},\boldsymbol{\ell},\mathbf{v})$ minimizing \eqref{eq:socp} corresponds to the underlying grid operation point, and more importantly, $f(\mathbf{p},\mathbf{q})$ is the actual power loss. Therefore, the instantaneous power loss $f_t(\mathbf{q}^g)=f(\mathbf{p}_t,\mathbf{q}^g-\mathbf{q}_t^c)$ has been expressed as the optimum value of an SOCP. Furthermore, since the function arguments $(\mathbf{p}_t,\mathbf{q}^g-\mathbf{q}_t^c)$ appear in the left-hand side of constraints \eqref{eq:sp}-\eqref{eq:sq}, $f_t(\mathbf{q}^g)$ is a perturbation function and is known to be convex~\cite[Lemma~4.24]{Ruszczynski}. 

The convexity of $f_t(\mathbf{q}^g)$ implies the existence of its subdifferential $\partial f_t(\mathbf{q}_t^g)$~\cite{Ruszczynski}. To efficiently calculate a subgradient $\mathbf{g}_t\in\partial f_t(\mathbf{q}_t^g)$, let us first eliminate $(\mathbf{P},\mathbf{v})$ and constraints \eqref{eq:sp} and \eqref{eq:sv} from \eqref{eq:socp}. To that end, define $\mathbf{z}:=[\mathbf{Q}'~\boldsymbol{\ell}']'$, and express $(\mathbf{P},\mathbf{v},\mathbf{q})$ as affine functions of $\mathbf{z}$, namely
\begin{subequations} \label{eq:lin}
\begin{align}
\mathbf{P}&=\mathbf{A}_p \mathbf{z} + \mathbf{b}_p(\mathbf{p})\label{eq:linP}\\
\mathbf{v}&=\mathbf{A}_v \mathbf{z} + \mathbf{b}_v\label{eq:linv}\\
\mathbf{q}&=\mathbf{A}_q \mathbf{z}\label{eq:linq}
\end{align}
\end{subequations}
for appropriate $N\times 2N$ matrices $\mathbf{A}_p,\mathbf{A}_v,\mathbf{A}_q$, and $N\times 1$ vectors $\mathbf{b}_p(\mathbf{p}), \mathbf{b}_v$. Notice the dependence of $\mathbf{b}_p$ on active injections $\mathbf{p}$. Using these substitutions, constraints \eqref{eq:sp} and \eqref{eq:sv} can be eliminated; the voltage constraints \eqref{eq:sV} can be expressed as $\underline{\mathbf{v}}\leq \mathbf{A}_v \mathbf{z} + \mathbf{b}_v \leq \overline{\mathbf{v}}$; and the equalities in \eqref{eq:sq} are compactly written as \eqref{eq:linq}. The $n$-th hyperbolic constraint in \eqref{eq:sl} can be expressed as the second-order cone~\cite{FCL}
\begin{equation}\label{eq:cone}
\left\|\left[\begin{array}{c}
2P_n\\
2Q_n\\
v_{\pi_n}-\ell_n
\end{array}
\right]\right\|_2\leq v_{\pi_n}+\ell_n
\end{equation}
or in terms of the introduced variable $\mathbf{z}$ as 
\begin{equation}\label{eq:cone2}
\|\mathbf{A}_n\mathbf{z}+\mathbf{b}_n(\mathbf{p})\|_2\leq \mathbf{c}_n'\mathbf{z}+d_n
\end{equation}
where the involved parameters are defined as
\begin{align*}
\mathbf{A}_n&:=\left[\begin{array}{c}
2\mathbf{e}_n'\mathbf{A}_p\\
2[\mathbf{e}_n'~\mathbf{0}']\\
\mathbf{e}_{\pi_n}'\mathbf{A}_v-[\mathbf{0}'~\mathbf{e}_n']
\end{array}
\right],&\mathbf{b}_n(\mathbf{p})&:=
\left[\begin{array}{c}
2\mathbf{e}_n'\mathbf{b}_p(\mathbf{p})\\
0\\
\mathbf{e}_{\pi_n}'\mathbf{b}_v
\end{array}\right]\nonumber\\
\mathbf{c}_n'&:=\mathbf{e}_{\pi_n}'\mathbf{A}_v+[\mathbf{0}'~\mathbf{e}_n'],& d_n&:=\mathbf{e}_{\pi_n}'\mathbf{b}_v.
\end{align*}
Using the aforementioned substitutions and for $\mathbf{r}_z:=[\mathbf{0}'~\mathbf{r}']'$ with $\mathbf{r}$ being the vector of line resistances, problem \eqref{eq:socp} can be equivalently written as
\begin{subequations}\label{eq:socp2}
\begin{align}
f(\mathbf{p},\mathbf{q})=\min_{\mathbf{z}}&~\mathbf{r}_z'\mathbf{z}\label{eq:s2cost}\\
\textrm{s.to}~&\mathbf{A}_q\mathbf{z}=\mathbf{q}\label{eq:s2q}\\
&\|\mathbf{A}_n\mathbf{z}+\mathbf{b}_n(\mathbf{p})\|_2\leq \mathbf{c}_n'\mathbf{z}+d_n,~n\in\mathcal{L}\label{eq:s2l}\\
&\underline{\mathbf{v}}\leq \mathbf{A}_v \mathbf{z} + \mathbf{b}_v \leq \overline{\mathbf{v}}\label{eq:s2v}
\end{align}
\end{subequations}
which is also an SOCP. Assumption (A2) and the fact that \eqref{eq:socp} is bounded below (by zero) guarantee strong duality and that the dual problem of \eqref{eq:socp2} is solvable~\cite[Proposition~5.3.2]{Be99}. Standard results from sensitivity analysis further imply that the subdifferential of $f(\mathbf{p},\mathbf{q})$ with respect to $\mathbf{q}$ coincides with the negative of the optimal dual variables corresponding to \eqref{eq:s2q}~\cite[Theorem 4.26]{Ruszczynski}, \cite{ConejoJOTA}.

The sought subgradient can be thus obtained via the dual problem of \eqref{eq:socp2}. Towards this direction, let $\boldsymbol{\lambda}$, $\underline{\boldsymbol{\nu}}\geq  \mathbf{0}$, and $\overline{\boldsymbol{\nu}}\geq  \mathbf{0}$, be the dual variables corresponding to \eqref{eq:s2q} and the lower and upper bounds in \eqref{eq:s2v}, respectively. To dualize the SOC constraints, introduce also the variable pairs $\{(\mathbf{u}_n,\mu_n)\}_{n\in\mathcal{L}}$. Then, the dual of \eqref{eq:socp2} is provided as~\cite{SOCP_Alizadeh},~\cite[pp.~566-7]{Be99}
\begin{align}
\max_{
\{\mathbf{u}_n,\mu_n\} \atop \boldsymbol{\lambda},\underline{\boldsymbol{\nu}},\overline{\boldsymbol{\nu}}}
&\sum_{n\in\mathcal{L}} \left(\mathbf{u}_n'\mathbf{b}_n(\mathbf{p}) -\mu_n d_n\right) -\boldsymbol{\lambda}'\mathbf{q} + \mathbf{b}_v'(\overline{\boldsymbol{\nu}} - \underline{\boldsymbol{\nu}})\label{eq:dual}\\
\textrm{s.to}~&\|\mathbf{u}_n\|_2\leq \mu_n\nonumber\\
&\underline{\boldsymbol{\nu}}\geq \mathbf{0},~\overline{\boldsymbol{\nu}}\geq \mathbf{0}\nonumber\\
&\mathbf{A}_v'(\overline{\boldsymbol{\nu}}-\underline{\boldsymbol{\nu}}) + \mathbf{A}_q'\boldsymbol{\lambda} + \sum_{n\in\mathcal{L}} \mathbf{A}_n'\mathbf{u}_n -\mu_n \mathbf{c}_n+\mathbf{r}_z=\mathbf{0}\nonumber
\end{align}
which can be solved as an SOCP as well. In deriving \eqref{eq:dual}, constraints \eqref{eq:s2l} have been dualized based on the fact that for fixed $(\mathbf{x},x_0)$, the maximization $\max_{\mathbf{u},\mu_0}\left\{\mathbf{u}'\mathbf{x}-\mu_0 x_0:~\|\mathbf{u}\|_2\leq \mu_0\right\}$ is equivalent to $\max_{\mu_0\geq 0} \mu_0(\|\mathbf{x}\|_2-x_0)$ and becomes zero when $\|\mathbf{x}\|_2 \leq x_0$; and infinity, otherwise.

If the tuple $\left(\{\mathbf{u}_{n,t}^{\star},\mu_{n,t}^{\star}\}_{n\in\mathcal{L}}, \boldsymbol{\lambda}_t^{\star},\underline{\boldsymbol{\nu}}_t^{\star},  \overline{\boldsymbol{\nu}}_t^{\star}\right)$ is a maximizer of \eqref{eq:dual} for $\mathbf{p}=\mathbf{p}_{t}$ and $\mathbf{q} = \hat{\mathbf{q}}_{t-1}^g-\mathbf{\mathbf{q}}_{t}^c$, then $-\boldsymbol{\lambda}_{t}^{\star}\in\partial f_t(\hat{\mathbf{q}}_{t-1}^g)$. Hence, the $\mathbf{g}_t$ in \eqref{eq:sa2} can be set to $\mathbf{g}_t=-\boldsymbol{\lambda}_t^{\star}$. Finally, under (A2), complementary slackness asserts that if $\mu_{n,t}^{\star} > 0$ for all $n\in\mathcal{L}$ maximizing \eqref{eq:dual}, then the related primal constraints in \eqref{eq:s2l} are satisfied with equality; see e.g.,~\cite{SOCP_Lobo}. Thus, $\{\mu_{n,t}^{\star}>0\}_{n\in\mathcal{L}}$ provides an exactness certificate for the convex relaxation in \eqref{eq:socp}.

Table~\ref{tbl:srpca} summarizes the novel stochastic reactive power compensation scheme, for which while two observations are in order.

\begin{remark}\label{re:merit1}
The derived control scheme does not depend on any distributional assumption on actual active and reactive power injections. It rather utilizes real-time microgrid operation data to infer the underlying statistics. The numerical tests in Section~\ref{sec:simulations} indicate that this data-driven approach can even track slow time-varying statistics.
\end{remark}

\begin{remark}\label{re:merit2}
Albeit the focus has been on minimizing the reactive power compensation cost, other microgrid management tasks (voltage deviation and conservation voltage regulation) could be amenable to this stochastic control framework.
\end{remark}


\begin{table}[t]
\renewcommand{\arraystretch}{1.05}
\caption{Stochastic Reactive Power Management Algorithm} \label{tbl:srpca}
\small
\begin{center}
\vspace{-1em}
\begin{tabular}{|p{0.9\linewidth}|}
\hline
1:\hspace*{0.5em} Input $\{c_n\}_{n\in\mathcal{N}}$, $(\underline{\mathbf{v}},\overline{\mathbf{v}})$, and step size $\eta_t>0$\\
2:\hspace*{0.5em} Construct $\mathbf{A}_p,\mathbf{A}_v,\mathbf{A}_q,\mathbf{b}_v,\mathbf{r}_z$, and $\{\mathbf{A}_n,\mathbf{c}_n,d_n\}_{n\in\mathcal{L}}$\\
3:\hspace*{0.5em} Initialize $\hat{\mathbf{q}}_{0}^g=\mathbf{0}$\\
4:\hspace*{0.5em} \textbf{for} $t=1,\ldots,T$ \textbf{do}\\
5:\hspace*{2em} Acquire $(\mathbf{p}_t,\mathbf{q}_t^c)$ and construct $\mathbf{b}_p(\mathbf{p}_t)$\\
6:\hspace*{2em} Solve \eqref{eq:dual} for $\mathbf{q}=\hat{\mathbf{q}}_{t-1}^g-\mathbf{q}_t^c$ to acquire $\boldsymbol{\lambda}_t^{\star}$\\
7:\hspace*{2em} Define $\mathbf{y}_t:=\hat{\mathbf{q}}_{t-1}^g + \eta_t\boldsymbol{\lambda}_t^{\star}$\\
8:\hspace*{2em} Apply $\hat{\mathbf{q}}_t^g$ as updated from \eqref{eq:closed}\\
9:\hspace*{0.5em} \textbf{end for}\\
\hline
\end{tabular}
\vspace{-1em}
\end{center}
\end{table}%

\subsection{Algorithm Convergence}\label{subsec:convergence}
Define the cost function in \eqref{eq:market} as
\begin{equation}\label{eq:marketcost}
h(\mathbf{q}^g):=\mathbb{E}[f_t(\mathbf{q}^g)] + \sum_{n=1}^N c_n |q_n^g|.
\end{equation}
The following result that can be obtained from \cite[Theorem~8]{COMID} characterizes the convergence of the iterates in \eqref{eq:sa}.

\begin{proposition}\label{pro:rate}
Let $\hat{\mathbf{q}}^g$ be a minimizer of \eqref{eq:market}, $\hat{\mathbf{q}}_t^g$ the update of \eqref{eq:sa}, and $\boldsymbol{\lambda}_t^{\star}$ a maximizer of \eqref{eq:dual}. If $\|\hat{\mathbf{q}}^g - \hat{\mathbf{q}}_t^g\|_2^2\leq 2 D^2$ and $\|\boldsymbol{\lambda}_t^{\star}\|_2\leq L$ for all $t$, it holds that 
\begin{equation}\label{eq:mean}
\mathbb{E}[h(\bar{\mathbf{q}}^g_T)] - h(\hat{\mathbf{q}}^g)\leq \frac{\alpha DL}{\sqrt{T}}
\end{equation}
where $\bar{\mathbf{q}}^g_T:=\tfrac{1}{T}\sum_{t=1}^T \hat{\mathbf{q}}_t^g$; and the constant $\alpha$ is 2 for $\eta_t=\frac{D}{L\sqrt{t}}$, and $3/2$ for $\eta_t=\frac{D}{L\sqrt{T}}$. It further holds that
\begin{equation}\label{eq:prob}
h(\bar{\mathbf{q}}^g_T) - h(\hat{\mathbf{q}}^g)  \leq 
\frac{DL}{\sqrt{T}}\left(\alpha + 4\sqrt{\log \delta}\right).
\end{equation}
with probability at least $1-\delta^{-1}$.
\end{proposition}
Proposition~\ref{pro:rate} guarantees that the expected power loss experienced by $\bar{\mathbf{q}}^g_T$ converges to the optimum stochastic power loss at the rate of $\mathcal{O}(1/\sqrt{T})$. Beyond mean value convergence from \eqref{eq:mean}, the bound in \eqref{eq:prob} assures that  $h(\bar{\mathbf{q}}^g_T)$ remains close to the optimum $h(\hat{\mathbf{q}}^g)$ with high probability. According to the online convex optimization terminology, the algorithm in  Table~\ref{tbl:srpca} enjoys sublinear regret~\cite{COMID}. Moreover, Proposition~\ref{pro:rate} asserts that the novel control scheme can operate for a constant step size $\eta_t=\tfrac{D}{L\sqrt{T}}$, assuming of course that $T$ is known in advance. That could be the case, if the proposed reactive power management scheme is periodically reset due to a new real-time active power market dispatch. If on the other hand, $T$ is unknown, a time-decaying step size $\eta_t=\tfrac{D}{L\sqrt{t}}$ works as well with a slight degradation in performance. Both the step sizes and the obtained bounds in Proposition~\ref{pro:rate} depend on $D$ and $L$. Apparently, when the reactive injection region $\mathcal{Q}$ models box constraints, $D$ depends on the reactive power capabilities of installed PVs as $D\leq 2\sum_{n\in\mathcal{N}_q}(\overline{q}_n^g)^2$. Regarding $L$, the $\ell_2$-norms of the subgradients $\mathbf{g}_t$ can be upper bounded too when $\mathcal{Q}$ is compact~\cite{COMID}. Knowing precisely $L$ may be practically unrealistic. Interestingly enough though, if the step size is $\eta_t=\frac{\beta D}{L \sqrt{t}}$ for some $\beta>0$ rather than $\eta_t=\frac{D}{L \sqrt{t}}$, then \eqref{eq:mean} holds for $\alpha=\tfrac{3}{2}\max\{\beta,1/\beta\}$~\cite{Nemirovski09}.

\section{Numerical Tests}\label{sec:simulations}
The novel stochastic reactive power management scheme is numerically tested first on a 47-bus industrial distribution network from South California Edison that is depicted in Fig.~\ref{fig:47grid}~\cite{GLTL12}. For each operation interval, the microgrid controller collects injections from load buses, as well as active injections from DG buses. Reactive injections from DG units are determined: (i) by solving the deterministic control scheme of \eqref{eq:instant}, and (ii) via the stochastic control scheme of Table~\ref{tbl:srpca}. Performance is tested in terms of the \emph{reactive power management cost} that is the instantaneous counterpart of the cost in \eqref{eq:market_tilde} evaluated on the true rather than the observed $(\mathbf{p}_t,\mathbf{q}_t^c)$. Observe that if $\{\tilde{c}_n=0\}_{n\in\mathcal{N}_q}$, the reactive power management cost coincides with the power loss cost. The power loss price is set to $\tilde{c}_0=6.6\cent$/kWh, and reactive power support prices are $\tilde{c}_n=\tilde{c}_0/80=0.0825\cent$/kVar \& h for all $n\in\mathcal{N}_q$. It is worth mentioning that all SOCP relaxations were feasible and exact.

\begin{figure}
\centering
\includegraphics[scale=0.16]{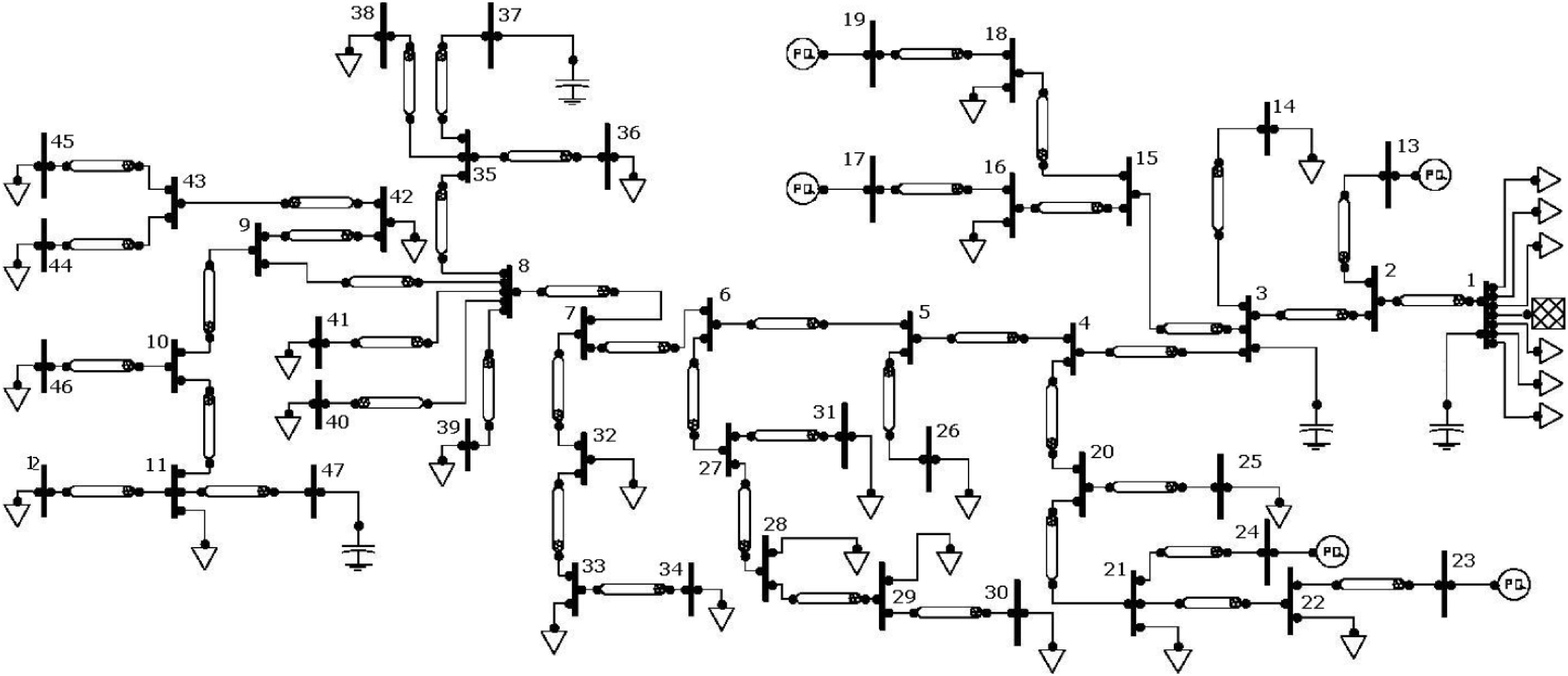}
\vspace*{-0.5em}
\caption{Schematic diagram of the 47-bus industrial distribution feeder with high penetration of photovoltaics located at buses 13, 17, 19, 23, and 24~\cite{GLTL12}.}
\vspace*{-0.5em}
\label{fig:47grid}
\end{figure}

\begin{figure}
\centering
\includegraphics[scale=0.60]{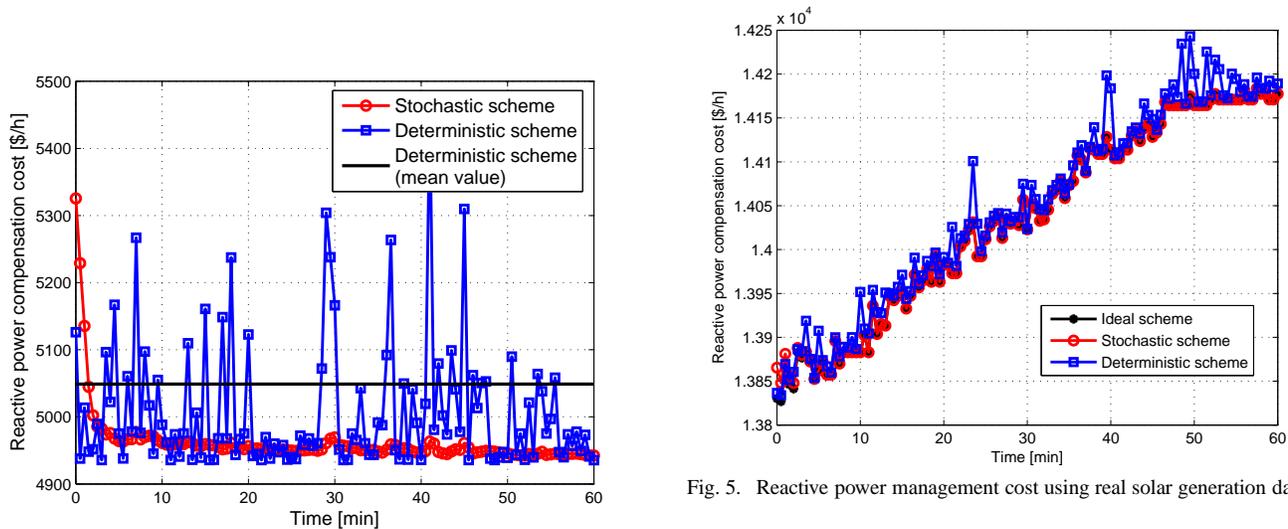}
\vspace*{-2em}
\caption{Reactive management cost for the 47-bus distribution grid $(\eta_t=2)$.}
\vspace*{-0.5em}
\label{fig:instantcost}
\end{figure}

The first experiment evaluates the effect of uncertainties in $(\mathbf{p}_t,\mathbf{q}_t^c)$ for controlling the 47-bus grid. Load injections $(\mathbf{p}^c,\mathbf{q}^c)$ are kept fixed throughout the interval to 45\% of their peak values with a power factor of 0.8. Photovoltaic injections $\mathbf{p}^g$ and shunt capacitors are kept fixed throughout the interval to 60\% of their peak values, while $v_0=1$. A period of 1 hour divided into 30-second control intervals is simulated. At each 30-sec interval, the controller observes a noise-corrupted version of the nominal $\mathbf{p}^o$ as $\mathbf{p}_{t}=\mathbf{p}^o+\boldsymbol{\epsilon}_t$, where the entries of $\boldsymbol{\epsilon}_t$ are independent and zero-mean Gaussian samples having variance $0.12$, thus modeling disturbances in power injections by 30\%. Noisy readings are likewise collected for the nominal $\mathbf{q}^{o,c}$. Although reactive PV injections $\hat{\mathbf{q}}_t^g$ are decided upon the noise-corrupted readings $(\mathbf{p}_t,\mathbf{q}_t^c)$, the actual power loss depends on their nominal values as $f(\mathbf{p}^o,\hat{\mathbf{q}}_t^g-\mathbf{q}^{o,c})$. The algorithm was implemented using MATLAB and CVX, and every reactive control was run within 1.2 secs on an Intel CPU @ 3.4 GHz (32 GB RAM) computer. Figure~\ref{fig:instantcost} depicts the reactive power management cost for the two control schemes over a single system realization. The algorithm of Table~\ref{tbl:srpca} converges within 20 iterations to a low cost, while its deterministic alternative fluctuates at consistently higher costs. 

\begin{figure}
\centering
\includegraphics[scale=0.35]{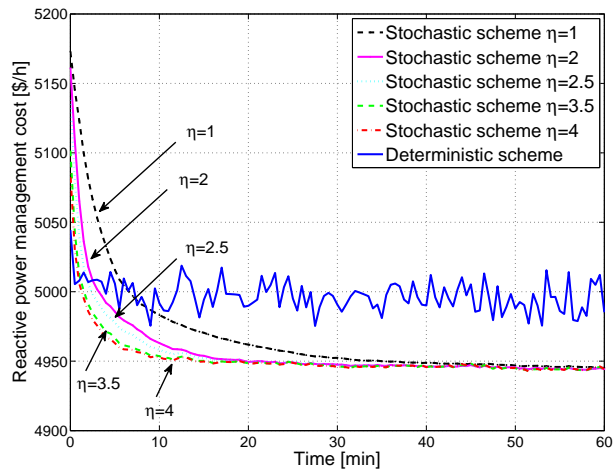}
\vspace*{-0.5em}
\caption{Reactive management cost averaged over 40 independent realizations of the 47-bus distribution feeder.}
\vspace*{-0.5em}
\label{fig:averagecost}
\end{figure}

Figure~\ref{fig:averagecost} presents the cost curves obtained after averaging 40 independent realizations. The curves verify that the stochastic scheme achieves significantly lower reactive power management costs than the myopic deterministic scheme. It is numerically observed that larger step sizes yield slower convergence, yet at a lower steady-state cost. The savings in \$/h are 28.7, 39.7, 41.8, 44.9, and 45.6, respectively, for $\eta=1,2,2.5,3.5,4$. Practically, tuning $\eta$ trades off the initial transient for the steady-state cost and the tracking of underlying statistics.

\begin{figure}
\centering
\includegraphics[scale=0.58]{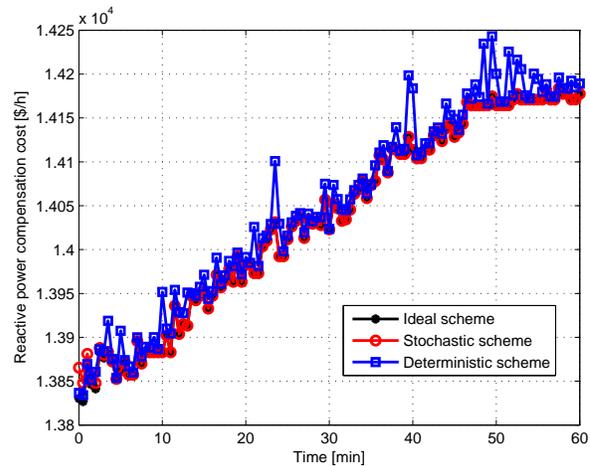}
\vspace*{-1.2em}
\caption{Reactive power management cost using real solar generation data~\cite{SmartStar}.}
\label{fig:realcost}
\vspace*{-1em}
\end{figure}

\begin{figure*}
\centering
\includegraphics[scale=0.38]{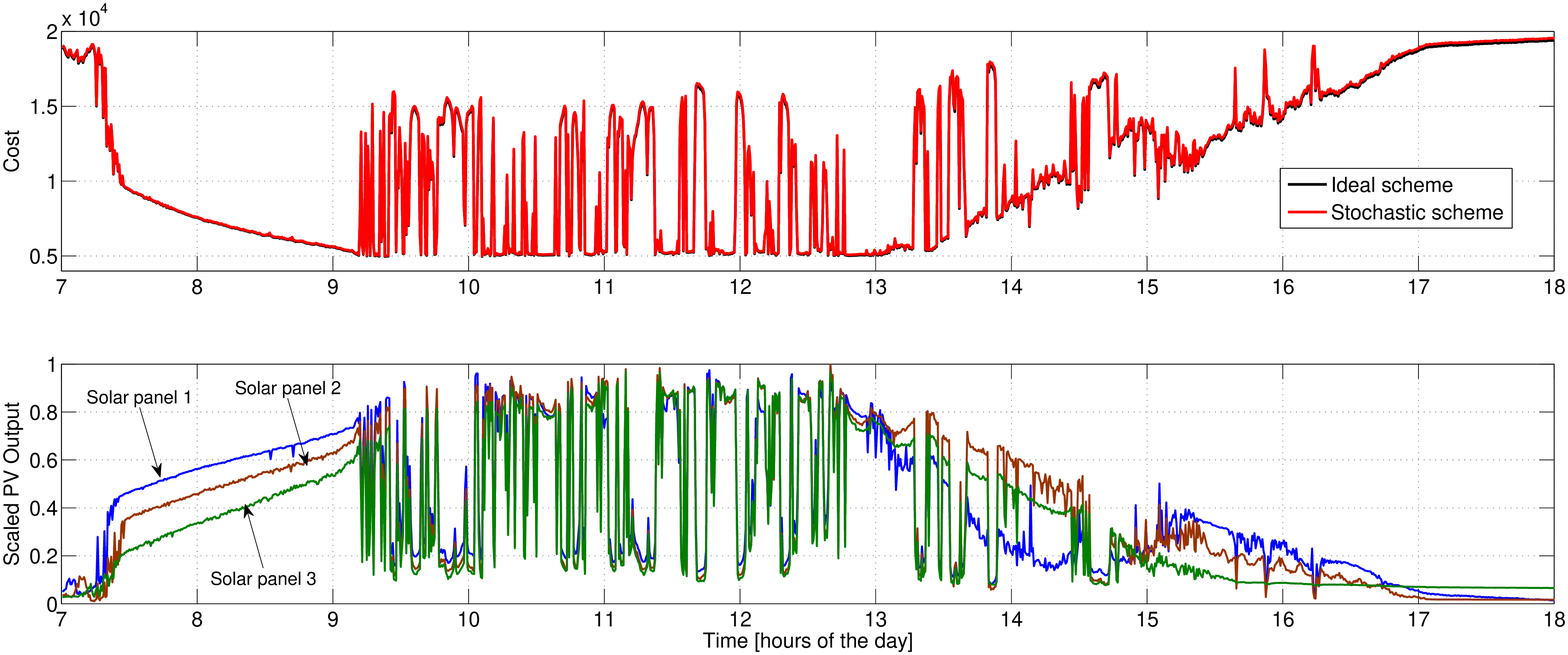}
\vspace*{-1em}
\caption{Top: Reactive power management cost with real solar generation data during August 12, 2011. Bottom: Normalized PV power output over 7am-6pm.}
\label{fig:movie}
\vspace*{-1em}
\end{figure*}

The second experiment entails real solar generation data from the Smart* project~\cite{SmartStar}. The power outputs of the 3 PVs involved in the Smart* microgrid over August 12, 2011, were preprocessed as follows: upon removing the minimum daily value, generation curves were normalized to the capacity of the PV units in Fig.~\ref{fig:47grid}; see also Fig.~\ref{fig:movie}. Industrial load demands were simulated at 80\% of their maximum values plus a Gaussian variation with standard deviation 15\% of the nominal value. In addition to the original PV generators on buses 13, 17, 19, 23, and 24; four more PV generators with capacity 1.2MW have been installed on buses 11, 28, 40, and 44, to model higher solar penetration. Figure~\ref{fig:realcost} shows the reactive power compensation cost attained over the period 18:30-19:30 at 30-sec control intervals. The controller determines the optimal reactive control based on the observed grid state which is the actual state delayed by 1 minute due to communication and computation delays. Together with the deterministic and stochastic schemes ($\eta_t=0.2$), the figure depicts the cost of the ideal control scheme that determines DG reactive injections based on the actual instantaneous grid state. Note that the latter is practically infeasible, but it serves as a lower bound. The numerical results show that upon convergence, the stochastic scheme approaches the ideal one and is able to track solar generation variations. The reactive power management benefit of the stochastic scheme over the deterministic one is 12.7\$/h.

Figure~\ref{fig:movie} presents the cost achieved by the new scheme during the daylight interval on August 12, 2011. Load demands were scaled to 90\% of their maximum value, and $\eta=0.2$. The control interval was selected as 30 seconds, and the observed state was the actual one delayed by 30 seconds. When PV generation is high, power losses and the related cost are low, as expected due to local generation. The curves on the top panel testify that the stochastic scheme attains a slightly higher cost than the ideal one. It further tracks successfully the steady solar power ramp occurring between 7.30-9.15am, as well as the variations due to cloud coverage for the rest of the day.

\begin{figure}
\centering
\includegraphics[scale=0.3]{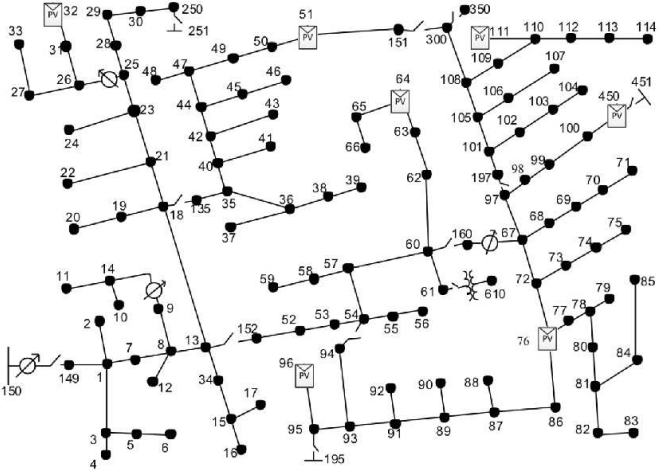}
\vspace*{-0.5em}
\caption{Schematic diagram of the IEEE 123-bus feeder with PVs~\cite{PSTCA}.}
\vspace*{-0.5em}
\label{fig:47grid}
\end{figure}

\begin{figure*}
\centering
\includegraphics[scale=0.38]{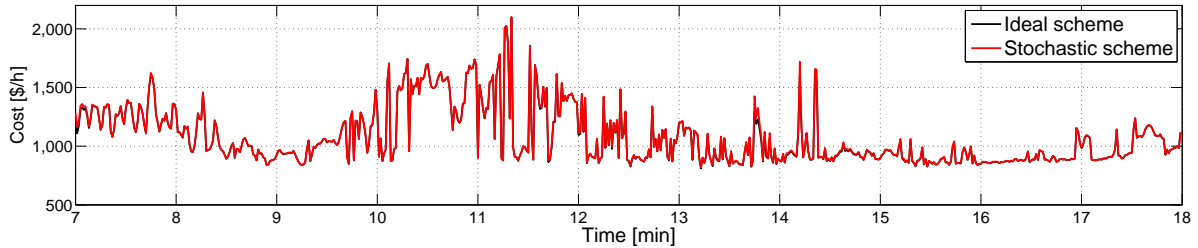}
\vspace*{-2em}
\caption{Reactive power management cost with real solar generation and demand data during August 12, 2011.}
\label{fig:movie2}
\vspace*{-1em}
\end{figure*}

Finally, the third experiment involved real data both for solar generation and consumption, which were tested on the IEEE 123-bus feeder~\cite{PSTCA}. The latter is a residential feeder that was modified according to~\cite{GLTL13}. Regarding renewable generation, solar panels were located on buses No. 32, 51, 64, 76, 96, 111, and 450, with capacities 40, 80, 160, 60, 160, 80, and 60~kW, respectively. Solar outputs were scaled versions of the curves shown at the bottom panel of Fig.~\ref{fig:movie}. All PV inverters were assumed to be oversized by 130\%, yielding a reactive power capacity of 0.66 times the active power capacity. Reactive power compensation prices were selected to be 0.0132 $\cent$/kVar \& h for all PVs. Consumption data provided by the Smart* project were utilized~\cite{SmartStar}. From a total of 443 minute-based household load data for August 12, 2011, we selected 85 after eliminating bad and incomplete entries. These load data were then scaled to match the IEEE 123-bus feeder profile. A reactive control period of one minute was implemented, and the observed states were equal to the actual ones delayed by one minute. The step size was set to $\eta_t=2$, and every control run lasted 3.7 secs. As shown in Fig.~\ref{fig:movie2}, the stochastic scheme was able to successfully track solar and load variations.

\section{Concluding Remarks}\label{sec:conclusions}
Reactive power compensation was considered in this work. Uncertainty and delays in acquiring microgrid states motivate well stochastic solutions. Building on a convex relaxation of the underlying problem as well as recent advances in online convex optimization, a novel stochastic scheme was developed. Reactive power injections from PV inverters were updated in real time. Numerical tests on practical microgrids verified that the novel control scheme converged within 10-20 iterations. The reactive power management cost attained was consistently lower than the one achieved by its myopic deterministic alternative. During experiments using real solar generation and load consumption data, the novel scheme tracked successfully the underlying system variations and approached the ideal reactive control scheme. The merit of our stochastic framework is twofold: First, apart from requiring slow variations, no distributional assumptions on active injections are imposed. Rather, the control algorithm adjusts dynamically to microgrid operation data. Second, albeit the goal here was to minimize the reactive power compensation cost, the novel approach could be extended to other pertinent microgrid management tasks (e.g., voltage deviation, conservation voltage regulation). Characterizing the resiliency of this centralized control to cyber-attacks on injection data and deriving decentralized solvers constitute directions for future research.

\appendix\label{sec:appendix}
\begin{IEEEproof}[Proof of Proposition~\ref{pro:closed}]
For notational simplicity, consider the canonical problem for $b>0$ and $\underline{x}< 0 < \overline{x}$
\begin{equation}\label{eq:canonical}
\hat{x}:=\arg\min\left\{\tfrac{1}{2}(x-a)^2 + b|x|:\underline{x}\leq x \leq \overline{x}\right\}.
\end{equation}
If $(\underline{\xi},\overline{\xi})$ are the optimal Lagrange multipliers for the box constraints in \eqref{eq:canonical}, the KKT conditions imply that there exists a subgradient $s(\hat{x})$ of $|x|$ at $\hat{x}$ satisfying
\begin{subequations}\label{eq:KKT}
\begin{align}
\hat{x}&=a-bs(\hat{x}) + \underline{\xi} -\overline{\xi}\label{eq:KKT1}\\
\underline{\xi}(\underline{x}-\hat{x})&=0\label{eq:KKT2}\\
\overline{\xi}(\hat{x}-\overline{x})&=0\label{eq:KKT3}\\
\underline{\xi},\overline{\xi}&\geq 0\label{eq:KKT4}\\
\underline{x}&\leq x \leq \overline{x}.\label{eq:KKT5}
\end{align}
\end{subequations}
For the subgradient of $|x|$, it holds that $s(x)=\sign(x)$ for $x\neq 0$, and $|s(x)|\leq 1$, otherwise. Depending on the sign of $\hat{x}$, three cases can be identified. If $\hat{x}>0>\underline{x}$, condition \eqref{eq:KKT2} yields $\underline{\xi}=0$ and \eqref{eq:KKT1} reads $\hat{x}=a-b-\overline{\xi}>0$. Two subcases can be now considered: Either $\hat{x}\in (0,\overline{x})$ implying that $\overline{\xi}=0$ and $\hat{x}=a-b$ when $a\in (b,b+\overline{x})$; or, $\hat{x}=\overline{x}$ implying $\overline{\xi}=a-b-\overline{x}$ when $a\geq b+\overline{x}$. The case of negative $\hat{x}$ can be treated similarly. In the third case where $\hat{x}=0$, conditions \eqref{eq:KKT1}-\eqref{eq:KKT3} yield $\underline{\xi}=\overline{\xi}=0$ and $a=bs(\hat{x})$. Since $|s(0)|\leq 1$, this third case occurs only when $a\in [-b,b]$.
\end{IEEEproof}

\bibliographystyle{IEEEtran}
\bibliography{myabrv,power}

\end{document}